# Consistent selection via the Lasso for high dimensional approximating regression models

### Florentina Bunea[*,1]

*Florida State University*

**Abstract:** In this article we investigate consistency of selection in regression models via the popular Lasso method. Here we depart from the traditional linear regression assumption and consider approximations of the regression function $f$ with elements of a given dictionary of $M$ functions. The target for consistency is the index set of those functions from this dictionary that realize the most parsimonious approximation to $f$ among all linear combinations belonging to an $L_2$ ball centered at $f$ and of radius $r_{n,M}^2$. In this framework we show that a consistent estimate of this index set can be derived via $\ell_1$ penalized least squares, with a data dependent penalty and with tuning sequence $r_{n,M} > \sqrt{\log(Mn)/n}$, where $n$ is the sample size. Our results hold for any $1 \leq M \leq n^\gamma$, for any $\gamma > 0$.

## Contents

## 1. Introduction

In this paper we show that the popular Lasso technique can be used for consistent feature selection in high dimensional approximating regression models. We consider the following framework. Given a random pair $(X, Y)$, we let $f(x) = E(Y|X = x)$ be the conditional mean function, henceforth called the regression function. We aim to reconstruct consistently a sparse approximation of $f$ via linear combinations of elements of a given dictionary of functions $\mathcal{F} = \{f_1, \ldots, f_M\}$. This reconstruction will be based on $(X_1, Y_1), \ldots, (X_n, Y_n)$, a sample of independent random pairs distributed as $(X, Y) \in (\mathcal{X}, \Re)$, where $\mathcal{X}$ is a Borel subset of $\Re^d$; all functions $f_j$ are defined on $\mathcal{X}$. Our aim expresses the belief that, in many instances, even if $M$ is

---

[*]Supported in part by NSF Grants DMS-04-06049 and DMS-07-06829.
[1]Department of Statistics, Florida State University, Tallahassee, Florida, USA, e-mail: flori@stat.fsu.edu

large, only a subset of $\mathcal{F}$ may be needed to approximate $f$ well. If that is the case, it may be of interest to determine whether this set can be estimated consistently via a computationally efficient method. The focus of this work is on consistent selection via the Lasso when the size of $\mathcal{F}$ grows polynomially with the sample size $n$, that is $M = n^\gamma$, for any $\gamma > 0$.

We begin by giving a number of examples of dictionaries $\mathcal{F}$ and associated consistency issues.

1. If $d = M$ and $f_j(X) = X_j$ for all $j$, one may be interested in identifying the subset of variables with linear combinations close to $f$. A familiar particular case is linear regression, where one assumes that $f(X) = \lambda'X$, with $\lambda \in \Re^M$ having non-zero components in positions corresponding to a set $J^* \subseteq \{1, \ldots, M\}$. Here we depart from this traditional equality assumption and consider the more realistic case where $f$ is not equal to, but can be well approximated by a linear combination of the given variables. We discuss this in detail in the next section.
2. Another problem of interest is that of finding consistently a sparse linear approximation of $f$ realized with elements from a large list of $M$ possibly competing estimators. These estimates may correspond to $M$ different methods of estimation, may be computed from $M$ different samples with the same mean function, or may correspond to $M$ different values of a tuning parameter of the same method. Instances of the latter arise in kernel based methods that require the choice of a grid of values for the bandwidth parameter or in Bayesian methods, where the specification of a grid of values for hyper-parameters is needed. A consistent identification of a subset of the estimates in these examples would validate the use of a particular restriction on an initially large grid. In such situations, when the elements of $\mathcal{F}$ are estimators, we will assume that they have been computed on samples independent of the one used for subset selection and treat them here as fixed functions.
3. A last example is the nonparametric estimation of $f$ from a collection of $M$ given basis functions, where only a subset may realize a good approximation of $f$, as described in the following subsection.

There exist a number of model selection methods that yield consistent subset selection in regression models. In discussing them a number of distinctions are needed.

The first one pertains to the evolution of the literature on model selection techniques in regression. One important cut-off point in this evolution seems to be the computational complexity of a particular method and, within this, the size of $M$ relative to $n$ plays a crucial role. If $M \leq n$, procedures based on various information criteria occupy an important place. They are referred to now as the BIC/AIC-type methods; we mention here the seminal works of ([1], [15]), the unifying theory of [2], and, various generalizations of these methods ([4], [7]). Such procedures can be easily implemented for small to moderate $M$. For larger values of $M$ multiple testing procedures, in particular of the FDR type (e.g., [3], [9]), or cross-validation with all its variants (holdout validation, $K$-fold) [21], are popular, but become more computationally complex as $M$ increases. If $M > n$ these techniques may become computationally intractable, unless they are used as part of a multiple-stage scheme. For a further overview on computational aspects in model selection, from a Bayesian perspective, see [11].

Whereas the above mentioned methods can still be used in very particular regression models when $M > n$, for instance, for sequence-space models, where model selection via BIC is equivalent to hard thresholding, they typically fail, computa-tionally, when $M$ is large. A standard solution in this case is to seek estimates that solve a certain class of convex optimization problems. Among the most popular estimates of this type in regression is the penalized least squares estimate with an $\ell_1$-type penalty (Lasso), which we describe in detail in the next section. In a Bayesian framework it can be derived from a Gaussian likelihood with a Laplace prior. Two important aspects set the $\ell_1$ regularized (Lasso) type estimators apart: they are easy and fast to compute; see [8], [13], [14], [18], among others, for efficient algorithms; and, if $M > n$, some components of the estimate will be set to zero, in finite samples, see, e.g., [13]. Therefore, via a one-step easily implementable procedure, one obtains subset selection even if $M > n$. To date, this method (or its variants) is the most widely used in regression problems of very high dimension, especially when dimension reduction is of interest.

The second distinction in discussing consistency of selection in regression is related to the target for consistency. Consistency of selection has been studied for all aforementioned techniques *only* in the following context, which we term parametric: the target for selection is typically an index set $J^*$ corresponding to the non-zero true regression coefficients, whereas the remaining coefficients are assumed to be *exactly* zero. An estimation method that uses the data and all $M$ elements $f_j$ to yield a subset $\hat{I}$ of indices such that $P(\hat{I} = J^*) \to 1$ for large $n$ is called a consistent method of selection.

In light of these two distinctions we give below a summary of the existing results on consistency of selection. They have all been established for the traditional parametric target $J^*$.

If $M \leq n$ and under appropriate assumptions all the above methods, or close variants, yield consistent subset selection for the parametric target $J^*$. References include those for AIC/BIC-type methods ([4], [10], and [22], among others), multiple testing procedures [5], cross-validation procedures [16], and Lasso-type procedures [24].

If $M > n$ consistency of selection has only been studied for Lasso-type estimators. Again, in the existing literature, the target is the standard target $J^*$. The results are limited. Meinshausen and Buhlmann [12] showed that $P(\hat{I} = J^*) \to 1$ in Gaussian graphical models, under assumptions that are tailored to models for which, in our notations, $(Y, X_1, \ldots, X_M) \sim N(0, \Sigma)$. Consistency of selection has been established when $M > n$, for fixed design linear regression models and a target set $J^*$ that corresponds to coefficients $\lambda_j^*$ that are assumed to be lower bounded by a sequence of order $O(n^{-\delta/2})$, for $0 < \delta < 1$ [23]. Similar results, under slightly different assumptions, have also been obtained for a three stage procedure [20]: in the first stage Lasso estimates are computed for a number of values of the tuning parameter, in the second step cross-validation is performed to select one Lasso estimate, and in the third one the model is refitted on the variables present in the selected Lasso estimate. We also refer to a related notion of consistency, in fixed design regression with Gaussian errors [19].

If $M > n$ consistent subset selection via the Lasso has not been investigated, to the best of our knowledge, in the general framework we describe in detail below. Within this framework, we extend the existing results to more general regression models on a random design and a more general target index set.

### 1.1. Beyond linear regression

Despite its practical appeal, the study of selection procedures that are consistent for target sets other than the classical one has received very little attention. Our target set will be defined relative to linear approximations of $f$ with elements of $\mathcal{F}$ with respect to the $L_2(\nu)$ norm $\|\cdot\|$, where we denote the probability measure of $X$ by $\nu$.

Formally, define

$$(1.1) \qquad \Lambda = \left\{ \lambda \in \Re^M : \, \|\sum_{j=1}^M \lambda_j f_j - f\|^2 \leq C_f r_{n,M}^2 \right\},$$

where $C_f > 0$ is a constant depending only on $f$ and $r_{n,M}$ is a positive sequence that converges to zero and which will be specified in the next section. In what follows we assume that $\Lambda$ is not void. For any $\lambda \in \Re^M$ we let $J(\lambda)$ denote the index set corresponding to the non-zero components of $\lambda$ and denote by $M(\lambda)$ its cardinality. Let $k^* = \min\{M(\lambda) : \lambda \in \Lambda\}$. We define our target vector

$$(1.2) \qquad \lambda^* = \operatorname{argmin}\left\{ \|\sum_{j=1}^M \lambda_j f_j - f\|^2 : \, \lambda \in \mathbb{R}^M, \, M(\lambda) = k^* \right\}.$$

Let $I^* = J(\lambda^*)$ denote the index set corresponding to the non-zero elements of $\lambda^*$ and note that $I^*$ has cardinality $k^*$. Thus $f^* = \sum_{j \in I^*} \lambda_j^* f_j$ provides the sparsest approximation to $f$ that can be realized with $\lambda \in \Lambda$ and, in particular,

$$(1.3) \qquad \|f^* - f\|^2 \leq C_f r_{n,M}^2.$$

This motivates our treating $I^*$ as the target index set.

We note that if one assumes, as in standard linear regression models, that $f(x) = \sum_{j=1}^M \lambda_j x_j = \sum_{j \in I^*} \lambda_j^* x_j = f^*(x)$, where $\lambda_j^*$ denotes the non-zero components of $\lambda$, then (1.3) is trivially satisfied for any positive sequence $r_{n,M}$. Therefore, the classical target $J^*$ is a particular case of ours.

In order to ensure that $\lambda^*$ captures the essential features of $f$ in a parsimonious way we require that its components not be unnecessarily small, otherwise we can place their indices outside $I^*$. Formally, we will require that the following condition holds.

**Condition (C).** There exists $B > 0$, independent of $n$ or $M$, such that

$$\min_{j \in I^*} |\lambda_j^*| > B r_{n,M}.$$

We show below that $\ell_1$ penalized least squares can be used to estimate consistently the new target $I^*$, even if $M$ is larger than $n$, in particular if it grows as $n^\gamma$, for any $\gamma > 0$, under minimal assumptions on the dictionary $\mathcal{F}$ and appropriate choices for $r_{n,M}$. In Section 2 below we introduce the estimate and discuss these choices. Section 2.1 contains our main result, Theorem 2.1, together with a discussion of the assumptions under which it holds. The proof of the main result is given in Section 2.2 and intermediate results are proved in the Appendix.

## 2. Consistent selection via $\ell_1$ penalized least squares

We estimate the set $I^*$ of the previous section via $\ell_1$ penalized least squares. We first compute

$$\widehat{\lambda} = \underset{\lambda \in \Re^M}{\arg\min} \left\{ \frac{1}{n}\sum_{i=1}^{n}\{Y_i - \sum_{j=1}^{M}\lambda_j f_j(X_i)\}^2 + \text{pen}(\lambda) \right\}, \qquad (2.4)$$

where

$$\text{pen}(\lambda) = 2\sum_{j=1}^{M}\omega_{n,j}|\lambda_j| \quad \text{with} \quad \omega_{n,j} = r_{n,M}\|f_j\|_n, \qquad (2.5)$$

for a sequence $r_{n,M}$ given below, where we write $\|g\|_n^2 = n^{-1}\sum_{i=1}^{n}g^2(X_i)$ for any function $g: \mathcal{X} \to \Re$. We note that each $\lambda_j$ in the penalty term has a different, data-dependent, weight. The estimate $\widehat{\lambda}$ thus obtained is in one-to-one correspondence with the following estimate. For each $1 \leq j \leq M$ define $\theta_j = 2\omega_{n,j}\lambda_j$ and let $A$ be the $M \times M$ diagonal matrix with diagonal entries $2\omega_{n,j}$. Next observe that $F\lambda = F_1\theta$, where $F$ is the $n \times M$ matrix with entries $f_j(X_i)$, $F_1 = FA^{-1}$ and $\theta = A\lambda$. Thus, denoting by $Y$ the $n$ dimensional vector with entries $Y_i$, the problem reduces to calculating

$$\widehat{\theta} = \underset{\theta \in \Re^M}{\arg\min} \frac{1}{n}(Y - F_1\theta)'(Y - F_1\theta) + \sum_{j=1}^{M}|\theta_j|,$$

for which the aforementioned fast algorithms can be used. Then, we compute our sought solution $\hat{\lambda} = A^{-1}\hat{\theta}$.

We let $\hat{I}$ denote the index set corresponding to the non-zero components of $\hat{\lambda}$. We show in the next subsection that $P(\hat{I} = I^*) \to 1$ when $n \to \infty$. We begin by noticing that we always have

$$P(\hat{I} = I^*) \geq 1 - P(I^* \not\subseteq \hat{I}) - P(\hat{I} \not\subseteq I^*).$$

Therefore, proving that $\hat{I}$ is consistent reduces to showing that each of the probabilities in the right-hand side of the inequality above converge to zero. In what follows we motivate choices for the sequence $r_{n,M}$ that stem from sufficient conditions under which this convergence is achieved. The proofs are presented in the next section.

We begin by noticing that if $\hat{\lambda} \to \lambda^*$, with probability converging to one, then $I^* \not\subseteq \hat{I}$ with probability converging to zero. To see this, further note that if component-wise consistency of $\hat{\lambda}$ holds, we will estimate *all* non-zero elements of $\lambda^*$ by non-zero sequences, but we may also estimate some of its zero components by some small, but non-zero sequences. In light of this fact, a first set of restrictions on $r_{n,M}$ will be such that $\hat{\lambda}$ is close to $\lambda^*$, in the sense below. It follows immediately (by [5], Theorem 2.3; see the Appendix below for a full formulation) that, with high probability

$$r_{n,M}|\hat{\lambda} - \lambda^*|_1 \leq D\{\|f - f^*\|^2 + k^* r_{n,M}^2\},$$

for some positive constant $D$, and where $|a|_1 = \sum_{j=1}^{M}|a_j|$ denotes the $\ell_1$ norm of any vector in $\Re^M$. Next, notice that the optimal parametric rate of convergence for a component $\hat{\lambda}_j$ of $\hat{\lambda}$ is of order $1/\sqrt{n}$, and it can be achieved if we knew $I^*$ of cardinality $k^* < M$ in advance. However, this is not known, so the best we can do is mimic this behavior in our context. We can do this by choosing $r_{n,M}$ of order $1/\sqrt{n}$, where we recall that we have assumed that $\|f - f^*\|^2 \leq r_{n,M}^2$. Notice further that this choice is optimal for the rate of convergence of $\hat{\lambda}$, which is not the focus here. Indeed, more modest rates of convergence of $\hat{\lambda}$ can be considered when consistency of selection is of main importance. We discuss in detail two concrete choices, and defer a complete analysis for future work.

One can consider $r_{n,M} = A\sqrt{\log(Mn)/n}$, for an appropriately large constant $A > 0$. Notice that this choice differs from the one that yields the optimal rate only by logarithmic factors, which are needed to accommodate dictionaries with $M > n$. With this choice, the target set $I^*$ corresponds to linear combinations of the elements of $\mathcal{F}$ that belong to, up to logarithmic factors, a $1/\sqrt{n}$ neighborhood of $f$, with respect to the $L_2(\nu)$ norm. This provides only a slight departure from the standard linear model assumption and standard target index set $J^*$. It is therefore not surprising that, in this case, our tuning sequence $r_{n,M}$ is also comparable to the one considered in parametric models ([12], [23]), where a sequence of the order of $1/n^{1/2-\theta}$, $\theta \in (0, 1/2)$, is employed. We note that this choice is slightly conservative, and can be relaxed to $O(\sqrt{\log(Mn)/n})$ in our framework, and therefore, as a particular case, in theirs.

In order to accommodate consistent selection in a purely nonparametric framework we need to increase the size of $r_{n,M}$. For instance, if all $f_j$ are estimates of $f$, and $r_{n,M}$ is as before, the set $\Lambda$ defined in (1.1) may be empty, as non-parametric estimates of $f$ have typically slower rates than $1/\sqrt{n}$. We therefore consider target sets $I^*$ corresponding to $L_2(\nu)$ neighborhoods around $f$ of radius $r_{n,M}^2$, now with $r_{n,M} = O\left((\log(Mn)/n)^{1/4}\right)$. In this case, the set $\Lambda$ given in (1.1) above is not empty if at least one of the estimators $f_j$ has, up to logarithmic factors, a rate of the order $n^{-1/4}$, which is a modest rate to require. Of course, if $f_j(X) = X_j$, as in linear regression, this choice means that we may be content with a coarser approximation than before. However, note that this approximation has the benefit of being realized with a smaller number of variables and that this may increase the interpretability of that particular model and be a desirable property in practical situations.

The results presented below hold for either of these choice, in particular for any $r_{n,M} \geq A\sqrt{\log(Mn)/n}$, and we will therefore not distinguish between them.

### *2.1. Main result: consistent subset selection*

We begin by listing and commenting on the assumptions under which our result holds. The first assumption refers to the error terms $W_i = Y_i - f(X_i)$. We recall that $f(X) = E(Y|X)$.

**Assumption (A1).** *The random variables $X_1, \ldots, X_n$ are independent, identically distributed random variables with probability measure $\mu$. The random variables $W_i$ are independently distributed with*

$$E\{W_i \mid X_1, \ldots, X_n\} = 0$$

*and*

$$E\{\exp(|W_i|) \mid X_1, \ldots, X_n\} \leq b \text{ for some finite } b > 0 \text{ and } i = 1, \ldots, n.$$

We also impose mild conditions on $f$ and on the functions $f_j$. Let $\|g\|_\infty = \sup_{x \in \mathcal{X}} |g(x)|$ for any function $g$ on $\mathcal{X}$.

**Assumption (A2).**

(a) There exists $0 < L < \infty$ such that $\|f_j\|_\infty \leq L$ for all $1 \leq j \leq M$.
(b) There exists $c_0 > 0$ such that $\|f_j\| \geq c_0$ for all $1 \leq j \leq M$.
(c) There exists $L_0 < \infty$ such that $E[f_i^2(X) f_j^2(X)] \leq L_0$ for all $1 \leq i, j \leq M$.
(d) There exists $L_1 < \infty$ such that $\|f\|_\infty \leq L_1 < \infty$.
(e) There exists $L^* < \infty$ such that $\|f - f^*\|_\infty \leq L^*$.

**Remark 2.1.** We note that $(a)$ trivially implies $(c)$. However, as the implied bound may be too large, we opted for stating $(c)$ separately. Note also that $(a)$ and $(d)$ imply the following: for any fixed $\lambda \in \Re^M$, there exists a positive constant $L(\lambda)$, depending on $\lambda$, such that $\|f - \sum_{j=1}^M \lambda_j f_j\|_\infty = L(\lambda)$. Inspection of the proof of Theorem 2.1 below shows that we can allow $L^*$ to grow very slowly with $n$. However, for sake of clarity in presentation we opted for treating it as fixed.

**Assumption (A3).** Let
$$\rho_M(i,j) = \frac{<f_i, f_j>}{\|f_i\| \|f_j\|},$$
where $<f_i, f_j> = Ef_i(X)f_j(X)$ and $\|f_i\| = E^{1/2} f_i^2(X)$. Assume that
$$\max_{i \in I^*} \max_{j \neq i} |\rho_M(i,j)| \leq \frac{C}{k^*},$$
for some constant $C > 0$.

**Remark 2.2.** Following [6], $C = 1/45$ is an allowable choice. Other choices are possible, but improvement of constants is beyond the scope of this paper.

**Remark 2.3.** Assumption (A3) reflects the belief that the correlations between functions $f_j$ with $j \in I^*$ and functions $f_j$ with $j \notin I^*$ should be small. However, we allow the correlations outside $I^*$ to be arbitrary. We note that this assumption replaces the standard orthonormality assumption on the design matrix: it is given in terms of theoretical quantities and it can hold even if $M > n$. It can be checked in practice by replacing the theoretical correlations by sample correlations.

We denote by $G$ the event that the $n \times M$ matrix $F$ with entries $f_j(X_i)$ has full rank. To avoid additional technicalities, the results of this paper can be regarded as conditional on $G$. Otherwise, all the results can be re-derived by intersecting all the relevant events with $G$ and $G^c$, under the additional assumption that $P(G^c)$ is appropriately small.

We can now state our main result which we prove in the next subsection.

**Theorem 2.1.** If assumptions (A1)–(A3) and condition (C) hold, and $k^* r_{n,M} \to 0$ then $P(\hat{I} = I^*) \to 0$.

**Remark 2.4.** The convergence above holds either if $M$ is fixed and $n \to \infty$ or if both $M, n \to \infty$, if $r_{n,M} \geq A\sqrt{\log(Mn)/n}$ for an appropriately large constant $A$. Therefore we obtain consistency for both choices of $r_{n,M}$ discussed above. In our derivations we require that $M$ does not grow faster than a power of $n$.

**Remark 2.5.** The condition $r_{n,M} k^* \to 0$ imposes restrictions on the size of $k^*$. If $r_{n,M} = O(\sqrt{\log(Mn)/n})$ the theorem above shows that we can recover consistently subsets of size $k^* = O(\sqrt{n}/\log n)$, up to other logarithmic factors. The choice $r_{n,M} = O(\log(Mn)/n)^{1/4}$ corresponds to a coarser approximation of $f$ than before, and the restriction on the number of approximating functions is now $k^* = O(n^{1/4}/\log n)$.

### 2.2. Proof of Theorem 2.1

Recall that

$$P(\hat{I} = I^*) \geq 1 - P(I^* \not\subseteq \hat{I}) - P(\hat{I} \not\subseteq I^*).$$

Therefore, proving that $\hat{I}$ is consistent reduces to showing that each of the probabilities in the right hand side of the inequality above converge to zero. We present this in the following two propositions. We defer the proof of the intermediate results to the Appendix.

**Proposition 2.2.** *If assumptions (A1)–(A3) and condition (C) hold, and $r_{n,M}k^* \to 0$, then $P(I^* \not\subseteq \hat{I}) \to 0$ as $n \to \infty$, for any $r_{n,M} \geq A\sqrt{\log(Mn)/n}$, with $A > 0$ large enough.*

*Proof.* We follow the same reasoning as [4]. Let $c_n = \min_{k \in I^*} |\lambda_k^*|$ and recall that $c_n > Br_{n,M}$, by condition (C). Therefore

$$\begin{aligned} P(I^* \not\subseteq \hat{I}) &\leq P(j \notin \hat{I} \text{ for some } j \in I^*) \\ &\leq P(|\widehat{\lambda}_j - \lambda_j^*| = |\lambda_j^*|) \\ &\leq P(|\widehat{\lambda}_j - \lambda_j^*| > c_n) \to 0, \text{ as } n \to \infty \end{aligned}$$

where, in the second inequality, we used that $\widehat{\lambda}_j = 0$ for $j \notin \hat{I}$, by the definition of $\hat{I}$. The last inequality follows from Corollary 1 presented in the Appendix below. □

**Proposition 2.3.** *If assumptions (A1)–(A3) hold and $r_{n,M}k^* \to 0$, then $P(\hat{I} \not\subseteq I^*) \to 0$, as $n \to \infty$, for any $r_{n,M} \geq A\sqrt{\log(Mn)/n}$, with $A > 0$ large enough.*

*Proof.* Let

$$h(\mu) = \frac{1}{n}\sum_{i=1}^n \{Y_i - \sum_{j \in I^*} \mu_j f_j(X_i)\}^2 + 2r_{n,M}\sum_{j \in I^*}||f_j||_n|\mu_j|,$$

and define

(2.6) $$\tilde{\mu} = \arg\min_{\mu \in \Re^{k^*}} h(\mu).$$

Let

$$\mathcal{B} = \bigcap_{k \notin I^*}\left\{|\frac{2}{n}\sum_{i=1}^n[Y_i - \sum_{j \in I^*}\tilde{\mu}_j f_j(X_i)]f_k(X_i)| < 2r_{n,M}||f_k||_n\right\}.$$

Let $\tilde{\lambda} \in \Re^M$ be the vector that has the components of $\tilde{\mu}$ in positions corresponding to the index set $I^*$ and components equal to zero otherwise. Thus, by abuse of notation, $\tilde{\lambda} = (\tilde{\mu}, 0)$. From Lemma 3.4 in the Appendix it follows that, on the set $\mathcal{B}$, $\tilde{\lambda}$ is a solution of (2.4). Recall that $\hat{\lambda}$ is a solution of (2.4) by construction. Then, by arguments similar to those used in ([13], Theorems 3.1 and 3.2) regarding the closeness of two solutions it follows that, on the set $\mathcal{B}$, $\hat{\lambda}_k = 0$ for $k \in I^{*c}$. Therefore $\hat{I} \subseteq I^*$ on the set $\mathcal{B}$. Hence

$$P(\hat{I} \not\subseteq I^*) \leq P(\mathcal{B}^c)$$
$$= P\left(\bigcup_{k \in \{1,\ldots,M\} \setminus I^*} \left\{|\frac{2}{n}\sum_{i=1}^n [Y_i - \sum_{j \in I^*} \tilde{\mu}_j f_j(X_i)] f_k(X_i)| \geq 2r_{n,M}\|f_k\|_n\right\}\right)$$
$$\leq \sum_{k \in \{1,\ldots,M\} \setminus I^*} P\left(\left\{|\frac{2}{n}\sum_{i=1}^n [Y_i - \sum_{j \in I^*} \tilde{\mu}_j f_j(X_i)] f_k(X_i)| \geq 2r_{n,M}\|f_k\|_n\right\}\right).$$

Let $k \in \{1,\ldots,M\} \setminus I^*$ be fixed. Define the sets

$$E_1(k) = \left\{\frac{1}{n}|\sum_{i=1}^n W_i f_k(X_i)| < r_{n,M}\|f_k\|_n/2\right\},$$
$$E_2(k) = \left\{\|f_k\|_n^2 \geq \frac{1}{4}\|f\|^2\right\},$$
$$E_3(k) = \left\{|\frac{1}{n}\sum_{i=1}^n f_j(X_i) f_k(X_i)| \leq 2|\langle f_j, f_k\rangle| + \delta_{n,M},\ j \in I^*\right\},$$

where $\delta_{n,M} = 2CL^2 r_{n,M}$ will be specified below. The choice of $\delta_{n,M}$ is purely technical and does not affect the overall results.

Let $\tilde{f} = \sum_{j \in I^*} \tilde{\mu}_j f_j$. Recall that $\lambda^* \in R^M$ given by (1.2) has zero components in positions corresponding to indices in $I^{*c}$, by definition. Let $\mu^*$ be the vector in $\Re^{k^*}$ obtained from $\lambda^*$ by deleting these zeros. Therefore $f^* = \sum_{j=1}^M \lambda_j^* f_j = \sum_{j \in I^*} \mu_j^* f_j$. By successive applications of the triangle inequality and since $\|f_k\|_n \leq L$, for all $k \in I^{*c}$, by assumption (A2) (a), we obtain:

$$(2.7) \quad P\left(\frac{1}{n}|\sum_{i=1}^n [Y_i - \sum_{j \in I^*} \tilde{\mu}_j f_j(X_i)] f_k(X_i)| \geq r_{n,M}\|f_k\|_n\right)$$
$$\leq P\left(\frac{1}{n}|\sum_{i=1}^n W_i f_k(X_i)| \geq r_{n,M}\|f_k\|_n/2\right)$$
$$+ P\left(\frac{1}{n}|\sum_{i=1}^n (f(X_i) - \tilde{f}(X_i)) f_k(X_i)| \geq r_{n,M}\|f_k\|_n/2\right)$$
$$\leq P(E_1^c(k)) + P\left(\frac{1}{n}|\sum_{i=1}^n (f^*(X_i) - \tilde{f}(X_i)) f_k(X_i)| \geq r_{n,M}\|f_k\|_n/4\right)$$
$$+ P\left(\frac{1}{n}\sum_{i=1}^n |(f(X_i) - f^*(X_i))| \geq r_{n,M}\|f_k\|_n/4L\right)$$
$$\leq P(E_1^c(k))$$
$$+ P\left(|(\sum_{j \in I^*}(\tilde{\mu}_j - \mu_j^*)\frac{1}{n}\sum_{i=1}^n f_j(X_i)) f_k(X_i)| \geq r_{n,M}\|f_k\|_n/4\right)$$
$$+ P\left(\frac{1}{n}\sum_{i=1}^n |(f(X_i) - f^*(X_i))| \geq r_{n,M}\|f_k\|_n/4L\right).$$

To bound the second term in the last inequality above we first notice that on the set $E_3(k)$ and under assumptions (A2) (a) and (A3) we have

$$|\sum_{j \in I^*}(\tilde{\mu}_j - \mu_j^*)\frac{1}{n}\sum_{i=1}^n f_j(X_i))f_k(X_i)|$$
$$\leq 2\sum_{j \in I^*}|\tilde{\mu}_j - \mu_j^*||\langle f_j, f_k\rangle| + \delta_{n,M}\sum_{j \in I^*}|\tilde{\mu}_j - \mu_j^*|$$
$$\leq \frac{2CL^2}{k^*}|\tilde{\mu} - \mu^*|_1 + \delta_{n,M}|\tilde{\mu} - \mu^*|_1.$$

Therefore, on $E_2(k) \cap E_3(k)$, and under assumptions (A2), (a) and (b), and (A3) we have

$$P(|(\sum_{j \in I^*}(\tilde{\mu}_j - \mu_j^*)\frac{1}{n}\sum_{i=1}^n f_j(X_i))f_k(X_i)| \geq r_{n,M}\|f_k\|_n/4)$$
$$\leq P(|\tilde{\mu} - \mu^*|_1 \geq \frac{c_0}{32CL^2}k^*r_{n,M}) + P(|\tilde{\mu} - \mu^*|_1 \geq \frac{c_0}{16}r_{n,M}\delta_{n,M}^{-1})$$
(2.8) $$\leq 2P(|\tilde{\mu} - \mu^*|_1 \geq \frac{c_0}{32CL^2}k^*r_{n,M}),$$

for $n$ large enough, since the assumption $k^*r_{n,M} \to 0$ implies that $k^*r_{n,M} \leq 1$ for large $n$, and we recall that we defined $\delta_{n,M} = 2CL^2 r_{n,M}$.

Lastly, notice that on the set $E_2(k)$ and under assumption (A2) (b) and (e) the third term of the last inequality in display (2.7) can be bounded by

(2.9) $$P(\frac{1}{n}\sum_{i=1}^n |(f(X_i) - f^*(X_i))| \geq \frac{c_0}{8L}r_{n,M}).$$

To complete the proof we need to show that $P(E_1^c(k)), P(E_2^c(k))$ and $P(E_3^c(k))$ and the probabilities in (2.8) and (2.9), when summed over $k \in \{1., \ldots, M\} \setminus I^*$, converge to zero as $n \to \infty$. We show this in Lemma 3.5, Corollary 2 and Lemma 3.6, respectively, in the Appendix below. This completes the proof of this result. □

## Appendix

In order to show Proposition 2.2 and to bound (2.8) above we will use twice ([6], Theorem 2.3 page 177) and we begin by stating it here, for completeness. For any $\lambda \in \Re^M$ we let $J(\lambda)$ denote the index set corresponding to the non-zero components of $\lambda$ and denote by $M(\lambda)$ its cardinality. Let $\rho(\lambda) = \max_{i \in J(\lambda)} \max_{j \neq i} |\rho_M(i,j)|$. With $\Lambda$ given by (1.1) in Section 1.1, let $\Lambda_1 = \{\lambda \in \Lambda : \rho(\lambda) \leq C/M(\lambda)\}$.

**Theorem 2.3** ([6]). *Assume that (A1) and (A2) hold. Then the $\ell_1$ penalized least squares estimator $\hat{\lambda}$ given by (2.4) satisfies, for any $\lambda \in \Lambda_1$*

(3.10) $$P\left\{|\hat{\lambda} - \lambda|_1 \leq B_1 r_{n,M}M(\lambda)\right\} \geq 1 - \pi_{n,M}(\lambda),$$

*where*

$$\pi_{n,M}(\lambda) \leq 14M^2 \exp\left(-c_1 n \min\left\{\frac{r_{n,M}^2}{L_0}, \frac{r_{n,M}}{L^2}, \frac{1}{L_0 M^2(\lambda)}, \frac{1}{L^2 M(\lambda)}\right\}\right)$$
$$\exp\left(-c_2 \frac{M(\lambda)}{L^2(\lambda)} nr_{n,M}^2\right)$$

for some positive constants $c_1, c_2$ depending on $c_0, C_f$ and $b$ only, and a constant $B_1$ depending on $c_0$ and $C_f$.

Notice now that by (1.3) and under assumption (A3), $\lambda^* \in \Lambda_1$. We therefore have the following corollary.

**Corollary 1.** *Assume that (A1)–(A3) hold. Then*

$$P\left\{|\widehat{\lambda}_j - \lambda_j^*| > B_1 r_{n,M}\right\} \leq \pi^*,$$

*for all $1 \leq j \leq M$, where $\pi^* = \pi_{n,M}(\lambda^*)$.*

*Proof.* From ([6], Theorem 2.3) we obtain

$$1 - \pi^* \leq P\left\{|\widehat{\lambda} - \lambda^*|_1 \leq B_1 k^* r_{n,M}\right\} \leq P\left\{\min_{1 \leq j \leq M}|\widehat{\lambda}_j - \lambda_j^*| \leq B_1 r_{n,M}\right\}.$$

This immediately implies the result. □

**Remark 3.1.** Notice that $\pi^* \to 0$ as $n \to 0$ for any $r_{n,M} \geq A\sqrt{\log(Mn)/n}$, and for $B = B_1$, as needed in Proposition 2.2 in Section 2.2 above.

In order to control the probability (2.8) we first define $U$ and $U_1$, the analogues of the sets $\Lambda$ and $\Lambda_1$ defined above.

$$U = \left\{\mu \in \Re^{k^*} : \|f - \sum_{j \in I^*} \mu_j f_j\|^2 \leq C_f r_{n,M}^2\right\}, \quad U_1 = \{\mu \in U : \rho(\mu)M(\mu) \leq C\}.$$

Recall that $\mu^*$ is the vector in $\Re^{k^*}$ obtained from $\lambda^*$ by deleting the zero entries. Then, since assumption (A3) implies $\max_{i \in I^*} \max_{j \in I^*, j \neq i} |\rho_M(i,j)| \leq C/k^*$ and $\|f - \sum_{j=1}^M \lambda_j^* f_j\| = \|f - \sum_{j \in I^*} \mu_j f_j\|$ we deduce that $\mu^* \in U_1$. Therefore, using again ([6], Theorem 2.3) applied now to the dictionary $\{f_j\}_{j \in I^*}$ and quantity $\widetilde{\mu}$ defined in (2.6) above, we obtain the following corollary:

**Corollary 2.** *Assume that (A1)–(A3) hold. Then*

(3.11) $$P\{|\widetilde{\mu} - \mu^*|_1 \leq B_2 k^* r_{n,M}\} \geq 1 - p^*,$$

*where*

$$p^* \leq 14 k^{*2} \exp\left(-c_1 n \min\left\{\frac{r_{n,M}^2}{L_0}, \frac{r_{n,M}}{L^2}, \frac{1}{L_0 k^{*2}}, \frac{1}{k^* L^2}\right\}\right)$$
$$+ \exp\left(-c_2 \frac{k^*}{L^2(\lambda^*)} n r_{n,M}^2\right),$$

*for some positive constants $c_1, c_2$ as above and a constant $B_2 > 0$ that only depends on $C_f$ and $c_0$.*

**Remark 3.2.** If $r_{n,M} \geq A\sqrt{\log(Mn)/n}$, then $Mp^* \to 0$ as $n \to \infty$, for $A > 0$ large enough. Hence, the probability given by (2.8), summed over $k$, converges to zero for both choices of $r_{n,M}$ introduced in Section 2, adjusting the value of $B_2$ if needed.

The following lemma is needed in the beginning of the proof of Proposition 2.3.

**Lemma 3.4.** $\tilde{\lambda} = (\tilde{\mu}, 0)$ *is a solution of* (2.4) *on the set*

$$\mathcal{B} = \bigcap_{k \notin I^*} \left\{ \left| \frac{2}{n} \sum_{i=1}^n [Y_i - \sum_{j \in I^*} \tilde{\mu}_j f_j(X_i)] f_k(X_i) \right| < 2r_{n,M} ||f_k||_n \right\}.$$

*Proof.* We recall that for any convex function $g : \Re^M \to \Re$ the subdifferential of $g$ at a point $\lambda$ is the set $D_\lambda = \{w \in \Re^M : g(u) - g(\lambda) \geq \langle w, u - \lambda \rangle\}$. Let $g(\lambda) = \frac{1}{n} \sum_{i=1}^n \{Y_i - \sum_{j=1}^M \lambda_j f(X_i)\}^2 + \text{pen}(\lambda)$, where we recall that our penalty term is $\text{pen}(\lambda) = 2r_{n,M} \sum_{j=1}^M ||f_j||_n |\lambda_j|$. Then (e.g., [13]) we have

$$D_\lambda = \{w \in \Re^M : w = -\frac{2}{n} F'(Y - F\lambda) + 2r_{n,M} v\},$$

where $v \in \Re^M$ is such that

$$\begin{aligned} v_k &= ||f_k||_n, & \text{if } \lambda_k > 0 \\ v_k &= -||f_k||_n, & \text{if } \lambda_k < 0 \\ v_k &\in [-||f_k||_n, ||f_k||_n], & \text{if } \lambda_k = 0, \end{aligned}$$

and where we recall that $Y = (Y_1, \ldots, Y_n)$ and $F$ is the $n \times M$ matrix with elements $f_j(X_i)$. By standard results in convex analysis, $\bar{\lambda} \in \Re^M$ is a point of local minimum for a convex function $g$ if and only if $0 \in D_{\bar{\lambda}}$, where $0 \in \Re^M$. Therefore, $\bar{\lambda}$ minimizes our $g(\lambda)$ if and only if $0 \in D_{\bar{\lambda}}$ if and only if

$$\left| \left( \frac{2}{n} F'(Y - F\bar{\lambda}) \right)_k \right| = 2r_{n,M} |v_k| \text{ for all } k \in \{1, \ldots, M\},$$

where $(\cdot)_k$ above denotes the $k$-th component of the vector in paranthesis. Equivalently, $\bar{\lambda}$ minimizes $g(\lambda)$ if and only if, for all $1 \leq k \leq M$

$$(3.12) \quad \left| \frac{2}{n} \sum_{i=1}^n [Y_i - \sum_{j=1}^M \bar{\lambda}_j f_j(X_i)] f_k(X_i) \right| = 2r_{n,M} ||f_k||_n, \text{ if } \bar{\lambda}_k \neq 0,$$

$$\left| \frac{2}{n} \sum_{i=1}^n [Y_i - \sum_{j=1}^M \bar{\lambda}_j f_j(X_i)] f_k(X_i) \right| \leq 2r_{n,M} ||f_k||_n, \text{ if } \bar{\lambda}_k = 0.$$

In what follows we find conditions under which $\tilde{\lambda} = (\tilde{\mu}, 0)$, with $\tilde{\mu}$ given in (2.6) above, satisfies (3.12). First notice that, by definition, $\sum_{i=1}^n [Y_i - \sum_{j=1}^M \tilde{\lambda}_j f_j(X_i)] = \sum_{i=1}^n [Y_i - \sum_{j \in I^*} \tilde{\mu}_j f_j(X_i)]$. Since $\tilde{\mu}$ is a solution of (2.6) then, by the above standard results in convex analysis, applied now to the function $h(\lambda)$ defined in the proof of Proposition 2.3, the following hold

$$\left| \frac{2}{n} \sum_{i=1}^n [Y_i - \sum_{j \in I^*} \tilde{\mu}_j f_j(X_i)] f_k(X_i) \right| = 2r_{n,M} ||f_k||_n, \text{ if } \tilde{\lambda}_k = \tilde{\mu}_k \neq 0, \ k \in I^*,$$

$$\left| \frac{2}{n} \sum_{i=1}^n [Y_i - \sum_{j \in I^*} \tilde{\mu}_j f_j(X_i)] f_k(X_i) \right| \leq 2r_{n,M} ||f_k||_n, \text{ if } \tilde{\lambda}_k = \tilde{\mu}_k = 0, \ k \in I^*.$$

Notice now that on the set $\mathcal{B}$ we also have

$$\left| \frac{2}{n} \sum_{i=1}^{n} [Y_i - \sum_{j \in I^*} \tilde{\mu}_j f_j(X_i)] f_k(X_i) \right| \leq 2r_{n,M} ||f_k||_n, \text{ if } k \notin I^* \text{ (for which } \tilde{\mu}_k = 0).$$

The above displays show that $\tilde{\lambda}$ satisfies condition (3.12) and is therefore a solution of (2.4) on $\mathcal{B}$. □

**Remark 3.3.** The observation that constitutes the statement of the above lemma has also been made elsewhere [12] for a slightly different penalty term. We have included here a full derivation of it for completeness and clarity.

To complete the proof of Proposition 2.3 we will make repeated use of Bernstein's inequality, which we state here for completeness.

**Bernstein's inequality.** *Let $\zeta_1, \ldots, \zeta_n$ be independent random variables such that*

$$\frac{1}{n} \sum_{i=1}^{n} E|\zeta_i|^m \leq \frac{m!}{2} w^2 d^{m-2}$$

*for some positive constants $w$ and $d$ and for all integers $m \geq 2$. Then, for any $\varepsilon > 0$ we have*

$$(3.13) \qquad P\left\{ \sum_{i=1}^{n} (\zeta_i - E\zeta_i) \geq n\varepsilon \right\} \leq \exp\left( -\frac{n\varepsilon^2}{2(w^2 + d\varepsilon)} \right).$$

**Lemma 3.5.** *Let assumptions (A1) and (A2) hold. Then*

$$\sum_{k \in \{1,\ldots,M\} \setminus I^*} P(E_1^c(k)) \to 0, \quad \sum_{k \in \{1,\ldots,M\} \setminus I^*} P(E_2^c(k)) \to 0, \text{ and}$$

$$\sum_{k \in \{1,\ldots,M\} \setminus I^*} P(E_3^c(k)) \to 0, \text{ as } n \to \infty.$$

*Proof.* To show $\sum_{k \in \{1,\ldots,M\} \setminus I^*} P(E_1^c(k)) \to 0$ it is enough to show that $(I) = \sum_{k \in \{1,\ldots,M\} \setminus I^*} P(E_1^c(k) \cap E_2(k)) \to 0$ and that $(II) = \sum_{k \in \{1,\ldots,M\} \setminus I^*} P(E_2^c(k)) \to 0$. The proofs follow immediately from Bernstein's inequality and the union bound. They are the same as ([6], proofs of Lemmas 4 and 5, page 186). We include here the derived probability bounds, for completeness.

$$(I) \leq 2M^2 \exp\left( -\frac{nr_{n,M}^2}{16b} \right) + 2M^2 \exp\left( -\frac{nr_{n,M} c_0}{8\sqrt{2}L} \right) + 2M^2 \exp\left( -\frac{nc_0^2}{12L^2} \right),$$

and

$$(II) \leq M^2 \exp\left( -\frac{nc_0^2}{12L^2} \right).$$

To bound the last quantity in the statement of the Lemma notice first that

$$\begin{aligned}
P(E_3^c(k)) &\leq 2\sum_{j\in I^*} P\left(\frac{1}{n}\sum_{i=1}^n f_j(X_i)f_k(X_i) > 2|\langle f_j, f_k\rangle| + \delta_{n,M}\right) \\
&\leq 2\sum_{j\in I^*} \exp\left\{-\frac{n}{4L_0}(|\langle f_j, f_k\rangle| + \delta_{n,M})^2\right\} \\
&\quad + 2\sum_{j\in I^*} \exp\left\{-\frac{n}{4L}(|\langle f_j, f_k\rangle| + \delta_{n,M})\right\} \\
&\leq 2M\exp\left\{-\frac{n\delta_{n,M}^2}{4L_0}\right\} + 2M\exp\left\{-\frac{n\delta_{n,M}}{4L}\right\}.
\end{aligned}$$

The second inequality of the display above follows from Bernstein's inequality with $\zeta_i = f_j(X_i)f_k(X_i)$, for every fixed $j$, and $k$ and with $w^2 = L_0$, $d = L^2$, for $\epsilon = |\langle f_j, f_k\rangle| + \delta_{n,m}$, used together with the inequality $e^{x/a+b} \leq e^{x/2a} + e^{x/2b}$ for all $x, a$ and $b$. Therefore, for $\delta_{n,M} = 2CL^2 r_{n,M}$ we obtain

$$\begin{aligned}
(III) &= \sum_{k\in\{1,\ldots,M\}\setminus I^*} P(E_2^3(k)) \\
&\leq 2M^2 \exp\left\{-\frac{C^2L^4 n r_{n,M}^2}{L_0}\right\} + 2M^2 \exp\left\{-\frac{CLnr_{n,M}}{2}\right\}.
\end{aligned}$$

Thus, the quantities $(I), (II)$ and $(III)$ converge to zero for any $r_{n,M} \geq A\sqrt{\log(M)n/n}$. □

**Lemma 3.6.** *Let assumptions (A1) and (A2) hold. Then*

$$(IV) = \sum_{k\in\{1,\ldots,M\}\setminus I^*} P\left(\frac{1}{n}\sum_{i=1}^n |(f(X_i) - f^*(X_i))| \geq \frac{c_0}{8L} r_{n,M}\right) \to 0.$$

*Proof.* By the Cauchy-Schwartz inequality we have

$$\begin{aligned}
(3.14)\quad & P\left(\frac{1}{n}\sum_{i=1}^n |(f(X_i) - f^*(X_i))| \geq \frac{c_0}{8L} r_{n,M}\right) \\
&\leq P\left(\frac{1}{n}\sum_{i=1}^n (f(X_i) - f^*(X_i))^2 \geq \frac{c_0^2}{64L^2} r_{n,M}^2\right) \\
(3.15)\quad &\leq P\left(\sum_{i=1}^n \{(f(X_i) - f^*(X_i))^2 - \|f - f^*\|^2\} \right.\\
&\qquad\left. \geq n\left(\frac{c_0^2}{64L^2} r_{n,M}^2 - \|f - f^*\|^2\right)\right) \\
&\leq P\left(\sum_{i=1}^n \{(f(X_i) - f^*(X_i))^2 - \|f - f^*\|^2\} \geq C_1 n r_{n,M}^2\right),
\end{aligned}$$

where we recall that $\|f - f^*\|^2 \leq C_f r_{n,M}^2$, by definition and $C_1 = c_0^2/64L^2 - C_f$, where we assume that we have already adjusted $C_f$ to have $C_1 > 0$, by taking an appropriate constant $A$ in the definition of $r_{n,M}$, if needed. The proof follows immediately from Bernstein's inequality applied to $\zeta_i = (f(X_i) - f^*(X_i))^2$, with $w = \sqrt{C_f} r_{n,M}$ and $d = L^*$, and for $\epsilon = C_1 r_{n,M}^2$. Therefore

$$(IV) \leq M \exp\{-\frac{C_f C_1^2}{4} n r_{n,M}^2\} + M \exp\{-\frac{C_1}{4L^*} n r_{n,M}^2\},$$

and both terms converge to zero for either choice of $r_{n,M}$. □